  \documentclass[a4paper,12pt]{amsart}
\usepackage{ae,amsfonts,euscript,enumerate}
\usepackage{amsmath,euscript}
\usepackage{amssymb}
\usepackage{amsthm}
\usepackage{enumerate}
\usepackage{amsfonts}
\usepackage{comment}
\usepackage{colortbl}
\usepackage{a4wide}
\usepackage{amssymb,amsmath,array}

\newtheorem{thm}{Theorem}[section]
\newtheorem{lem}[thm]{Lemma}

\theoremstyle{defn}

\def\co{{\mathcal O}}

\def\oqmm13{\co_q(M_{1,3})}
\def\oqm23{\co_q(M_{2,3})}

\input{xy}
\xyoption{all}

\def\eqref#1{(\ref{#1})}
\def\qed{~\vrule height8pt width 5pt depth -1pt\medskip}

\usepackage[dvips]{graphicx}
\usepackage{amssymb}
\usepackage{amsmath}
\usepackage[latin1]{inputenc}
\newif\ifpdf
       \ifx\pdfoutput\undefined
       \pdffalse 
       \else
       \pdfoutput=1 
       \pdftrue
\fi

\ifpdf
       \else
       \usepackage{graphicx,color,rotating}
\fi
\usepackage{fullpage}
\title[Centralizers in domains]{Centralizers in domains of finite Gelfand-Kirillov dimension}
\author{Jason P.~Bell}
\thanks{I thank NSERC for its generous support.  }

\address{Jason Bell\\
Department of Mathematics\\
Simon Fraser University\\
Burnaby, BC V5A 1S6
}

\email{jpb@math.sfu.ca}

\begin{document}
\bibliographystyle{plain}

\begin{abstract}
We study centralizers of elements in domains.  We generalize a result of the author and Small \cite{BS}, showing that if $A$ is a finitely generated noetherian domain and $a\in A$ is not algebraic over the extended centre of $A$ then the centralizer of $a$ has Gelfand-Kirillov dimension at most one less than the Gelfand-Kirillov dimension of $A$.  In the case that $A$ is a finitely generated noetherian domain of GK dimension $3$ over the complex numbers, we show that the centralizer of an element $a\in A$ that is not algebraic over the extended centre of $A$ satisfies a polynomial identity. 
\end{abstract}
\maketitle

\section{Introduction} 
We consider domains of finite Gelfand-Kirillov dimension. This dimension, defined by Gelfand and Kirillov, is defined as follows. 
Given a field $K$ and a finitely generated $K$-algebra $A$, we define the Gelfand-Kirillov dimension (GK dimension) 
of $A$ to be 
\[{\rm GKdim}(A) \ := \ \limsup_{n\rightarrow\infty} \frac{\log ({\rm dim}\,V^n)}{ \log\, n},\] 
where $V$ is a finite-dimensional subspace of $A$ that generates $A$ as a $K$-algebra. 
We note that this definition is independent of choice of $V$.

 In the case that $A$ is 
not finitely generated over $K$, we define the GK dimension of $A$ to be the supremum of 
the GK dimensions of all finitely generated subalgebras of A. Basic facts about 
GK dimension can be found in the excellent book of Krause and Lenagan \cite{KL}. 

Gelfand-Kirillov dimension is a noncommutative analogue of Krull dimension, and for this reason it has seen great application in the development of noncommutative algebraic geometry.  Just as the classification of curves and the classification of birational classification of surfaces have been of tremendous importance in algebraic geometry, so too have been the corresponding classifications thus far obtained for algebras of low GK dimension. 

A finitely generated domain of GK dimension $0$ is a division ring that is finite-dimensional over its center; A finitely generated domain of GK dimension $1$ over a 
field $K$ is a finite module over its center, and hence satisfies a polynomial identity \cite{SW}; if, 
in addition, $K$ is algebraically closed, a routine application of Tsen's theorem shows 
that this domain is in fact commutative.  Domains of GK dimension 2 are not well-understood in general, but Artin and Stafford have given a concrete description of these objects when the additional restriction that these domains be finitely generated $\mathbb{N}$-graded algebras that are generated in degree $1$.   

We consider the centralizers of elements of 
a finitely generated domain $A$ of finite GK dimension.   In the case that $A$ has GK dimension at most $2$, the centralizers are understood---they are all PI \cite{BS, Smok}.   The study of centralizers in fact has a much longer history, due to its connection to the theory of differential equations.  In particular, centralizers in the first Weyl 
algebra---a domain of GK dimension $2$---were studied as long ago 
as $1922$ by Burchnall and Chaundy \cite{BC}; ultimately, Amitsur showed that centralizers of non-scalar elements are isomorphic to polynomial algebras in one variable, in the case that the base field is the complex numbers \cite{Am}.

Our aim is to view these results under the lens of GK dimension; all results obtained thus far show that under reasonable hypotheses the GK dimension of the centralizer of an element in a domain has smaller GK dimension than the domain itself.  In the case that the domain has GK dimension $2$, this means that the GK dimension of the centralizer must be $1$ under reasonable hypotheses and so by well-known results in low GK dimension, the centralizer must be commutative in the case that the base field is algebraically closed.  We extend this result as follows.
To state our extension, we must introduce a few ideas that are used.  A semiprime noetherian algebra $A$ has a semisimple artinian quotient ring $Q(A)$; in the case that $A$ is a domain, $Q(A)$ is a division algebra.  We let $Z(Q(A))$ denote the centre of this ring (often called the \emph{extended centre} of $A$).  In general, if $a\in A$ is algebraic over $Z(Q(A))$ then we expect the centralizer of $a$ to be very large.  If, on the other hand, $a$ is not algebraic over $Z(Q(A))$, then we can show that the centralizer has smaller GK dimension.  

Given a ring $R$ and $a\in R$, we let $C(a;R)$ denote the centralizer of $a$ in $R$.  Our main result is given by the following theorem.

\begin{thm} Let $K$ be a field and let $A$ be a finitely generated $K$-algebra that is a domain of finite GK dimension.  If $a\in A$ is not algebraic over $Z(Q(A))$, the centre of the quotient division algebra of $A$, then
$${\rm GKdim}(C(a;A))\le {\rm GKdim}(A)-1.$$\label{thm: main1}
\end{thm}
As a corollary we prove the following result, which is a $3$-dimension analogue of a result of the author and Small \cite{BS} that says in a finitely generated domain of GK dimension $2$ that does not satisfy a polynomial identity, the centralizers of non-scalar elements are all commutative.
\begin{thm} Let $A$ be a finitely generated noetherian $\mathbb{C}$-algebra that is a domain of GK dimension $3$.  If $a\in A$ is not algebraic over $Z(Q(A))$, then $C(a;A)$ satisfies a polynomial identity.\label{thm: main2}
\end{thm}
\section{Proofs}
In this section, we give the proofs of our main results.
\begin{lem}  (Lenagan---appearing in Letzter \cite[Theorem 2.1]{Letzter}) .Let $S$ and $R$ be noetherian rings and suppose that $S$ is finite and free as a left and right $R$-module.  If $S$ is simple, then $R$ is simple. \label{lem: 1}
\end{lem}
\noindent {\bf Proof.} Suppose that $S$ is simple and $R$ is not simple.  Then $R$ has a nonzero proper two-sided ideal $I$.  Since $S$ is free as a left $R$-module, we can lift $I$ to a proper right ideal $J$ of $S$.  Then $M:=S/J$ is an $(R,S)$-bimodule.  Since $S$ is finite as a left $R$-module, $M$ is finitely generated as a left $R$-module.  Since $M$ is also a torsion right $S$-module, $M$ has a nonzero annihilator $I'$, which is a two-sided ideal of $S$, a contradiction. \qed
\vskip 2mm
\noindent {\bf Proof of Theorem \ref{thm: main1}.}  We may replace with the algebra $A$ with $Z(Q(A))A\subseteq Q(A)$ if necessary and so it is no loss of generality to assume that the base field $K=Z(Q(A))$.
 We note that $A$ is an ore domain since it has finite GK dimension.  Let $B=C(a;A)$.  If ${\rm GKdim}(B)>{\rm GKdim}(A)-1$, then by a theorem of Borho and Kraft \cite{BK}, $Q(A)$ is finite dimensional as both a left and right $Q(B)$-vector space.  
Let $R=Q(B)\otimes_K K(t)$ and $S=Q(A)\otimes_K K(t)$.  We note that $S$ is simple and is finite and free as a left $R$-module; moreover, $S$ and $R$ are both noetherian by the Hilbert basis theorem as they are localizations of polynomial extensions of division algebras.  By hypothesis, $a$ is not algebraic over $K$ and $a\in Z(Q(B))$.  Hence the centre of $Q(B)$ is an extension of $K$ of transcendence degree at least $1$.  In particular, the centre of $R$ contains a proper two-sided ideal generated by $a\otimes 1 -1\otimes t$ and hence is not a field; thus $R$ is not simple, contradicting Lemma \ref{lem: 1}.   \qed
\vskip 2mm
\noindent {\bf Proof of Theorem \ref{thm: main2}.} Let $B=C(a;A)$.  By Theorem \ref{thm: main1}, $B$ has GK dimension at most $2$.  If the GK dimension is less than $2$ then $B$ has GK dimension $1$ by Bergman's gap theorem as it contains the polynomial algebra $\mathbb{C}[a]$ \cite{KL}.   A complex domain of GK dimension $1$ is commutative by Tsen's theorem and the Small and Warfield theorem \cite{SW}, and so we may assume that $B$ has GK dimension $2$. 

We note that any finitely generated subalgebra $B'$ of $B$ is PI by a result of Smith and Zhang \cite{SZ}, as it has GK dimension at most $2$ and the centre contains $\mathbb{C}[a]$, which has GK dimension $1$.  In particular, $Q(B)$ is locally PI and has GK dimension $1$ as a $\mathbb{C}(a)$-algebra.

We note that $A$ is remains noetherian under extension of scalars as it is finitely generated over the complex numbers \cite{Bell}; hence $Q(A)$ also has this property, as $Q(A)\otimes_{\mathbb{C}} K$ is a localization of $A\otimes_{\mathbb{C}} K$ for any extension field $K$ of $\mathbb{C}$.  Since $Q(A)\otimes_{\mathbb{C}} K$ is free over $Q(B)\otimes_{\mathbb{C}} K$, we see that $Q(B)$ also remains noetherian under extension of scalars.  Let $K$ be an algebraically closed extension of $\mathbb{C}(a)$ with ${\rm Card}(K)>{\rm dim}_{\mathbb{C}(a)}(Q(B))$.  Then $R:=Q(B)\otimes_{\mathbb C(a)} K$ is a factor of the noetherian algebra $Q(B)\otimes_{\mathbb{C}} K$ and hence $R$ is noetherian.  Furthermore, any finitely generated subalgebra of $R$ is contained in $Q(B')\otimes_{\mathbb C(a)} K$ for some finitely generated subalgebra $B'$ of $B$; hence $R$ is locally PI and has GK dimension $1$ as a $K$-algebra.

Since the cardinality of $K$ is greater than ${\rm dim}_K(R)$, we see that the Jacobson radical of $R$ is nil; since $R$ is noetherian, it is nilpotent.  Let $N$ denote the Jacobson radical of $R$.

By Goldie's theorem, $$Q(R/N)\cong \prod_{i=1}^d M_{n_i}(D_i),$$ for some division algebras $D_1,\ldots ,D_d$ and natural numbers $n_1,\ldots , n_d$.  We note that each $D_i$ is locally PI and has GK dimension at most $1$ as a $K$-algebra as $R$ has this property.  Since $K$ is algebraically closed, each $D_i$ is commutative by Tsen's theorem and the Small-Warfield theorem \cite{SW}.   Thus $Q(R/N)$ is PI.  Since $N$ is nilpotent, $R$ is also PI.  Since $B$ embeds in $R$, we see that $B$ is PI.  The result follows. \qed 

\section{Concluding remarks}
In this section we make a few remarks about the results in this paper.

We note that if $A=\mathbb{C}[x_1,\ldots ,x_d][t,t^{-1};\sigma]$, where $\sigma(x_i)=\lambda_i x_i$ with $\lambda_1,\ldots ,\lambda_d$ generating a free abelian multiplicative subgroup of $\mathbb{C}^*$ of rank $d$, then $A$ is a noetherian domain of GK dimension $d+1$ and the centralizer of $x_1$ is the polynomial ring $\mathbb{C}[x_1,\ldots ,x_d]$, which has GK dimension $d$.  Thus it is possible for equality to occur in the conclusion of the statement of Theorem \ref{thm: main1}.  

We note that if $A$ is a noetherian domain of GK dimension $4$, the centralizer of an element that is not algebraic over $Z(Q(A))$ need not be PI.  To see this, let $A=\mathbb{C}[x_1,x_2][y_1,y_2;\delta_1,\delta_2]$ be the second Weyl algebra.  That is, we have relations given by $[x_i,y_j]=\delta_{i,j}$ and $[x_1,x_2]=[y_1,y_2]=0$.  Then the centralizer of $x_1$ is the subalgebra generated by $x_1,x_2$, and $y_2$; that is $C(x_1;A)$ is ismorphic to a Weyl algebra over $\mathbb{C}[x_1]$, which does not satisfy a polynomial identity.

We note that it is impossible to impose any bound on the PI degree of $C(a;A)$ in the conclusion of the statement of Theorem \ref{thm: main2}.  To see this, let 
$A=\mathbb{C}_q[x,y][t,t^{-1};\sigma]$, where $q$ is a primitive $n$'th root of unity and we have the relations $xy=qyx$, $\sigma(x)=2x$, $\sigma(y)=3y$.  Then $A$ is a noetherian $\mathbb{C}$-algebra of GK dimension $3$ and the centralizer of $x^n$ is $\mathbb{C}_q[x,y]$, whose quotient division algebra is $n^2$ dimensional over its centre, $\mathbb{C}(x,y)$.  Thus the centralizer satisfies a polynomial identity of degree $2n$ and does not satisfy any smaller identity by the Amitsur-Levitzki theorem.  

It would be nice to remove the hypotheses that the algebra be noetherian and that the base field be the complex numbers in the statement of Theorem \ref{thm: main2}---this would require a completely different approach than the one used here, which uses the fact that finitely generated noetherian algebras over the complex numbers are \emph{stably noetherian}; i.e., they remain noetherian under tensoring with field extensions of the base field.  In fact, the conclusion of the statement of Theorem \ref{thm: main2} holds over any base field if we insist that $A$ be \emph{stably noetherian} instead of simply being noetherian. 
\section*{Acknowledgments}
I thank Dan Rogalski, Toby Stafford, and Lance Small for many interesting conversations.

\end{document}